\documentclass[12pt, twoside]{article}
\usepackage{amsmath, amsfonts, amsthm, amssymb,epsfig}
\usepackage{epsfig}
\usepackage{graphicx}
\usepackage{color}
\usepackage{epstopdf}
\usepackage{fancyhdr}
\usepackage{amsthm}

\theoremstyle{definition}
\newtheorem*{theorem}{Theorem}
\newtheorem*{remark}{Remark}
\newtheorem*{definition}{Definition}
\newtheorem*{problem}{Problem}
\newtheorem*{example}{Example}
\newtheorem*{proposition}{Proposition}
\newtheorem*{lemma}{Lemma}
\newcommand\A{{\mathcal A}}

\newcommand{\p}{\partial}

\fancypagestyle{mypagestyle}{
  \fancyhf{}  
  \fancyhead[OC]{\textit{\small{Invariants of Links}}}
  \fancyhead[EL]{\thepage}
  \fancyhead[EC]{\small{J\'ozef H. PRZYTYCKI and Pawe{\l} TRACZYK}}
  
}
\fancypagestyle{firststyle}
{
   \fancyhf{}
    \lhead{PRZYTYCKI, J. H. and TRACZYK, P. \\ KOBE J. MATH,\\ 4(1987) 115-139}
    
}
\title{\centerline{\large \bf INVARIANTS OF LINKS OF CONWAY TYPE}}
\author{\normalsize By J\'ozef H. Przytycki and Pawe{\l} Traczyk} 
\date{\small(Received April 8, 1985)}
\begin{document}
\maketitle
\thispagestyle{firststyle}
\section{Introduction}The purpose of this paper is to present a certain combinatorial method of constructing invariants of isotopy classes of oriented tame links. 
This arises as a generalization of the known polynomial invariants of Conway and Jones. These invariants have one striking common feature. 
If $L_+, L_-$ and $L_0$ are diagrams of oriented links which are identical, except near one crossing point, 
where they look like in Fig. 1.1\footnote{Added for e-print: we follow here the old Conway's convention \cite{C};
in modern literature the role of $L_+$ and $L_-$ is usually inverted.  In \cite{Prz} the new convention is already used.}, 
then $w_{L_+}$ is uniquely determined by $w_{L_-}$ and $w_{L_0}$, and also $w_{L_-}$ is uniquely determined by $w_{L_+}$ and $w_{L_0}$. 
Here $w$ denotes the considered invariants (we will often take the liberty of speaking about the value of 
an invariant for a specific link or diagram of a link rather than for an isotopy class of links). In the above context, 
we agree to write $L_+^p, L_-^p$ and $L_0^p$ if we need the crossing point to be explicitly specified.
\ \\ \ \\
\centerline{{\psfig{figure=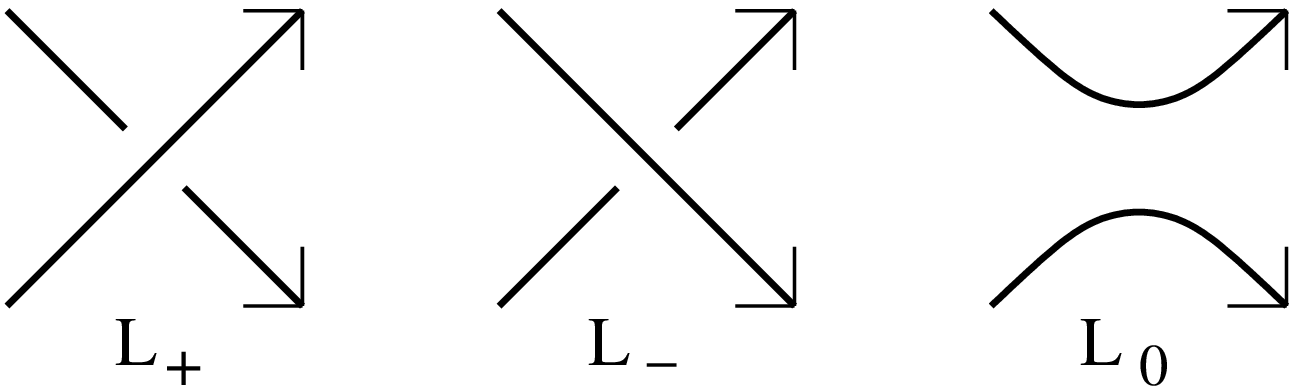,height=3.5cm}}}\ \\ 
\centerline{\footnotesize{Fig. 1.1}}
\indent In this paper we will consider the following general situation. Assume we are given an algebra $\A$ with a countable number of 
$0$-argument operations $a_1, a_2,..., a_n,...$ and two 2-argument operations $|$ and $*$. We would like to construct invariants
 satisfying the conditions 
\begin{align*}
    w_{L_+}= &\ w_{L_-} | w_{L_0} \text{ and}\\
    w_{L_-}= &\ w_{L_+} * w_{L_0} \text{ and} \\
    w_{T_n}= &\ a_n \text{ for $T_n$ being a trivial link of n components.}
\end{align*}
\indent We say that $(\A, a_1, a_2,..., |, *)$ is a Conway algebra if the following conditions 
are satisfied\footnote{Added for e-print: We were unaware, when writing this paper, that the condition 1.5,
$(a\!*\!b)\!*\!(c\!*\!d) =  (a\!*\!c)\!*\!(b\!*\!d)$, was used already for over 50 years, first time in \cite{Bu-Ma}, under various names, e.g.
entropic condition (see for example \cite{N-P}).}: \\
${1.1} \quad  a_n | a_{n+1} = a_n   \\
{1.2} \quad  a_n * a_{n+1} =  a_n  $\\
$\left.
\begin{aligned}
{1.3} &&  &(a|b)|(c|d) = (a|c)|(b|d) \\
{1.4} &&   &(a|b)\!*\!(c|d) =  (a\!*\!c)|(b\!*\!d) \\
{1.5} &&  &(a\!*\!b)\!*\!(c\!*\!d) =  (a\!*\!c)\!*\!(b\!*\!d)
\end{aligned}
\right\}
\quad \text{transposition properties}$ \\
${1.6} \quad   (a|b)*b = a \\ 
{1.7} \quad  (a*b)|b = a.$\\
\vspace{1mm} \\
We will prove the following theorem: \\
\begin{theorem}\textbf{1.8.} For a given Conway algebra $\A$ there exists a uniquely determined invariant $w$ which attaches an element 
$w_L$ from $\A$ to every isotopy class of oriented links and satisfies the conditions\\
${(1)} \quad w_{T_n} = a_n \hspace{2.5cm} \text{ initial conditions}$\\
$\left.
\begin{aligned}
{(2)}&& &w_{L_+}  = &\ w_{L_-} | w_{L_0}\\
{(3)}&&  &w_{L_-}  = &\ w_{L_+} * w_{L_0}
\end{aligned}
\right.\bigg\}
\ Conway\ relations$
\end{theorem}
It will be proved in \S{2}. \\
\indent Let us write here a few words about the geometrical meaning of the axioms 1.1-1.7 of Conway algebra. Relations 1.1 and 1.2 are introduced to 
reflect the following geometrical relations between the diagrams of trivial links of $n$ and $n+1$ components: 

\centerline{\epsfig{figure=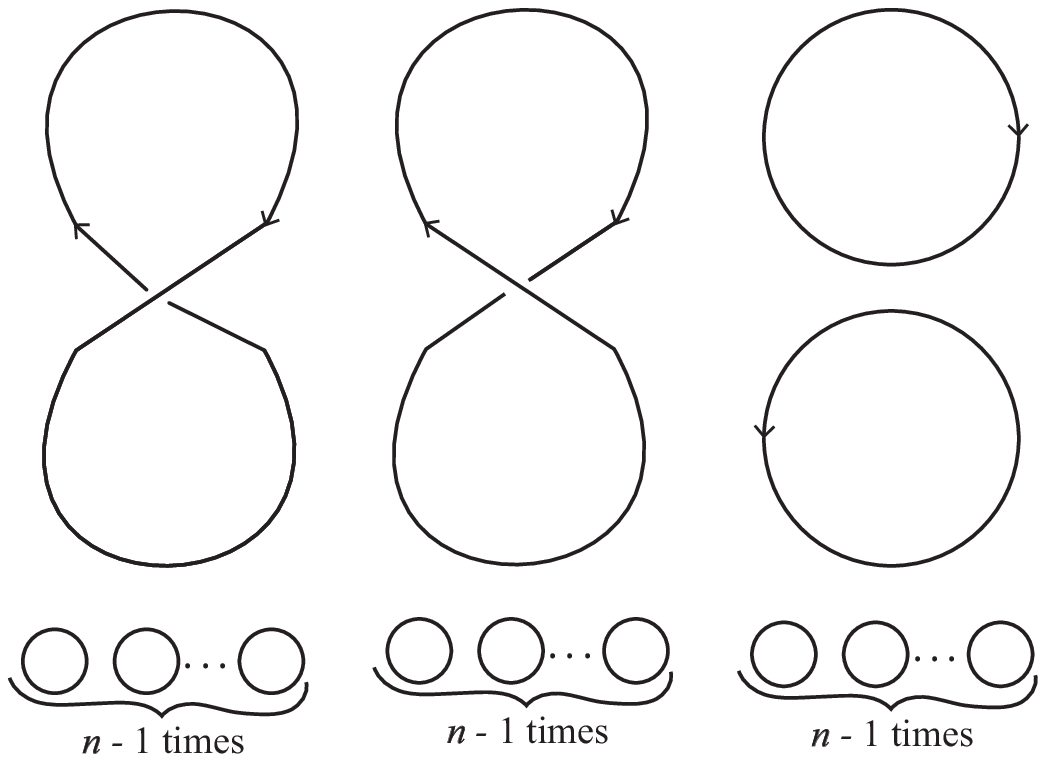,height=10cm}}
\centerline{\footnotesize{Fig. 1.2}}
Relations 1.3, 1.4, and 1.5 arise when considering rearranging a link at two crossings of the diagram but in different order. It will become clear in \S{2}. Relations 1.6 and 1.7 reflect the fact that we need the operations $|$ and $*$ to be in some sense opposite one to another. \\
\begin{example}\textbf{1.9.} (Number of components). Set $\A=N-$ the set of natural numbers; $a_i=i$ and $i|j=i*j=i$. 
This algebra yields the number of components of a link.\end{example}
\begin{example}\textbf{1.10.} Set $\A=$ \big\{0, 1, 2\big\}; the operation $*$ is equal to $|$ and 
$0|0\!=\!1,\ 1|0\!=\!0,\ 2|0\!=\!2,\ 0|1\!=\!0,\ 1|1\!=\!2,\ 2|1\!=\!1,\ 0|2\!=\!2,\ 1|2\!=\!1,\ 2|2\!=\!0$. Furthermore $a_i \equiv i \mod3.$ 
The invariant defined by this algebra distinguishes, for example, the trefoil knot from the trivial knot. \end{example}
\begin{example}\textbf{1.11.}{(a)} $\A=Z [x^{\mp 1} , y^{\mp 1}, z]$; $a_1 = 1, a_2=x+y+z,\dots, a_i = (x+y)^{i-1} + z(x+y)^{i-2} +
\dots+z(x+y)+z=(x+y)^{i-1}+z \big( \frac{(x+y)^{i-1}-1}{x+y-1} \big)$,{\small ...} . 
We define $|$ and $*$ as follows: $w_2 | w_0 = w_1$ and $w_1 * w_0 = w_2$ where
\begin{flalign*}
    {1.12}&& xw_1+yw_2 &=w_0-z, \ \  w_1,w_2,w_0 \in \A.&
\end{flalign*}
 \end{example}
\indent (b) $\A = Z$[ $x^{\mp1}$ , $y^{\mp1}$] is obtained from the algebra described in (a) by substitution $z=0$. In particular $a_i=(x+y)^{i-1}$ and 1.12 reduces to:
\begin{flalign*}
    {1.13}&& xw_1+yw_2 &= w_0.&
\end{flalign*}
We describe this algebra for two reasons: \\ 
\indent $-$first, the invariant of links defined by this algebra is the simplest generalization of the Conway polynomial 
\big(substitute $x=- \frac{1}{z}$, $y=\frac{1}{z}$, {[$K$-1]}\big) and the Jones polynomial 
\big(substitute $x= \frac{1}{t} \frac{1}{\sqrt{t}-\frac{1}{\sqrt{t}}}$, $y= \frac{-t}{\sqrt{t}-\frac{1}{\sqrt{t}}}$ {[$J$]}\big); \\
\indent$-$second, this invariant behaves well under disjoint and connected sums of links:
\begin{center}
    $\qquad P_{L_1 \sqcup L_2} (x,y)\!=\! (x+y)P_{L_1{\sharp} L_2}(x,y)\!=\!(x+y)P_{L_1}(x,y)\!\cdot\!P_{L_2}(x,y)$ \\
    \vspace{2mm}
    where $P_L\!(x,y)$ is a polynomial invariant of links yielded by $\A$.\\
    \vspace{2mm}
    \begin{example} \textbf{1.14.} (Linking number). Set $\A\!=\!N\!\times\!Z$, $a_i\!=\!(i,0)$ and\end{example}
\begin{equation*}
(a,b)|(c,d)=\begin{cases}
    (a,b-1) \ \textup{if} \ a > c \\
    (a,b) \qquad \textup{if} \ a \leqslant c 
\end{cases}
\end{equation*}
\begin{equation*}
(a,b)*(c,d)=\begin{cases}
    (a,b+1)\ \textup{if} \ a > c \\
    (a,b)\ \qquad \textup{if}\ a \leqslant c 
\end{cases}
\end{equation*}
\end{center}
\indent The invariant associated to a link is a pair (number of components, linking number).
\begin{remark}\textbf{1.15.} It may happen that for each pair $u,v\!\in\!\A$
there exists exactly one $w\!\in\!\A$ such that $v|w\!=\!u$ and $u\!*\!w\!=\!v$. Then we can introduce a new operation $\circ$: 
$\A \times \A \rightarrow \A$ putting $u \circ v=w$ (we have such a situation in Examples $1.10$ and $1.11$ 
but not in $1.9$ where $2|1\!=\!2\!*\!1\!=\!2\!=\!2|3\!=\!2\!*\!3$). Then $a_n\!=\!a_{n-1}\circ a_{n-1}$. If the operation $\circ$ 
is well defined we can find an easy formula for invariants of connected and disjoint sums of links. 
We can interpret $\circ$ as follows: if $w_1$ is the invariant of $L_+$ (Fig 1.1) and $w_2$ of $L_-$ then $w_1 \circ w_2$ 
is the invariant associated to $L_0$. \end{remark}
\begin{remark}\textbf{1.16.} Our invariants often allow us to distinguish between a link and its mirror image. 
If $P_L (x, y, z)$ is an invariant of $L$ from Example 1.11 (a) and $\overline{L}$ is the mirror image of $L$ then
\[P_{\overline{L}}(x,y,z) = P_L(y,x,z).\]
 \end{remark}
\indent We will call a crossing of the type \parbox{1.1cm}{\psfig{figure=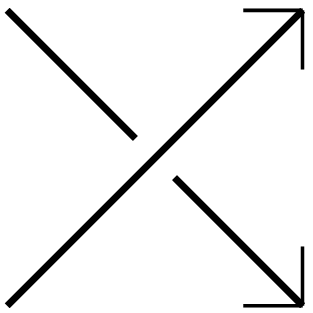,height=0.9cm}}
positive and crossing of the type 
\parbox{1.1cm}{\psfig{figure=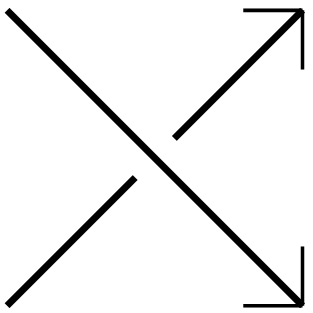,height=0.9cm}}
negative. This will be denoted by sgn $p= +$ or $-$. Let us consider now the following example.
\begin{example}\textbf{1.17.} Let $L$ be the figure eight knot represented by the diagram \end{example}
\centerline{\epsfig{figure=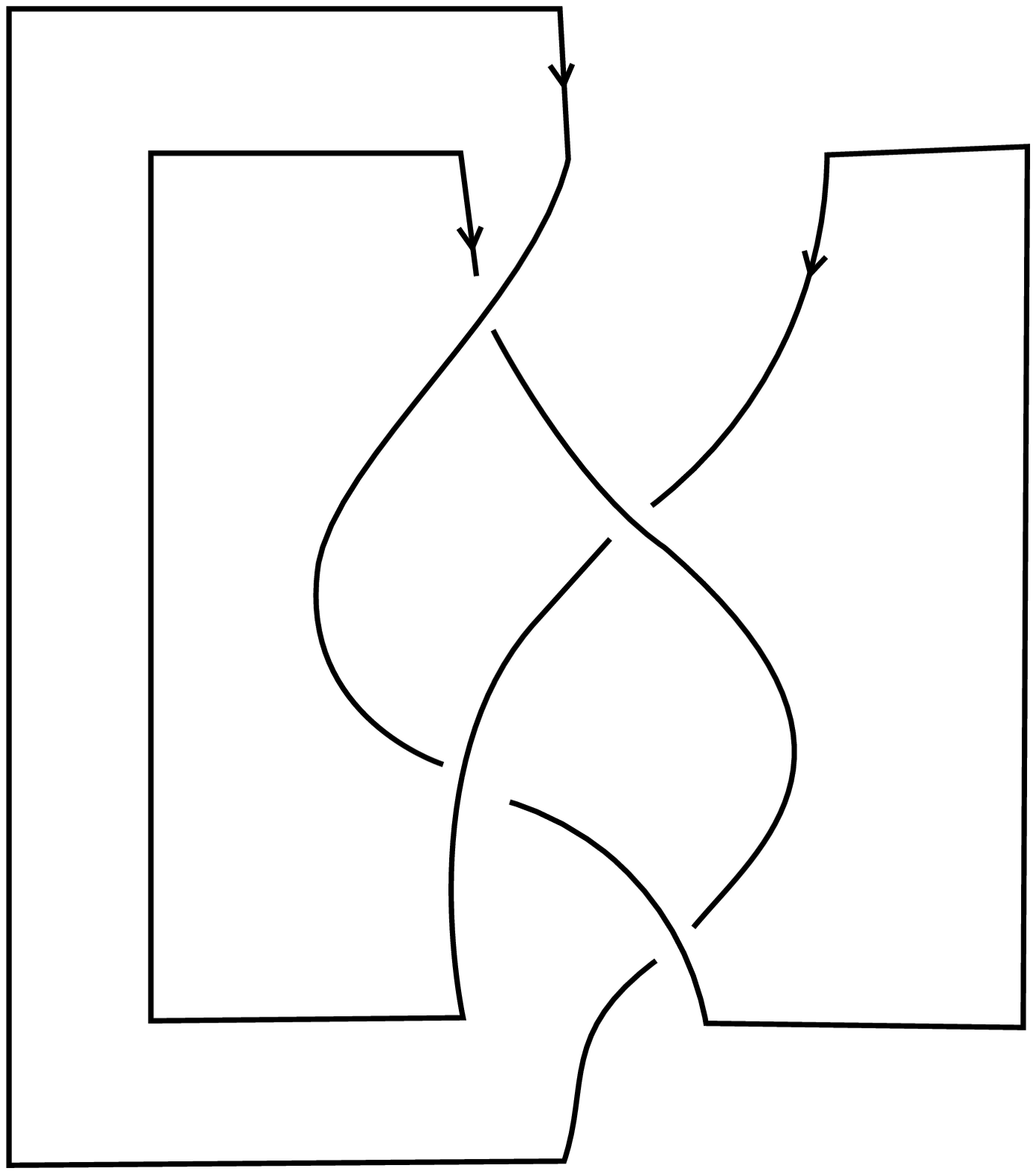,height=7.5cm}} \ \\ 
\centerline{\footnotesize{Fig. 1.3}} \ \\ 
To determine $w_L$ let us consider the following binary tree: \\
\centerline{\epsfig{figure=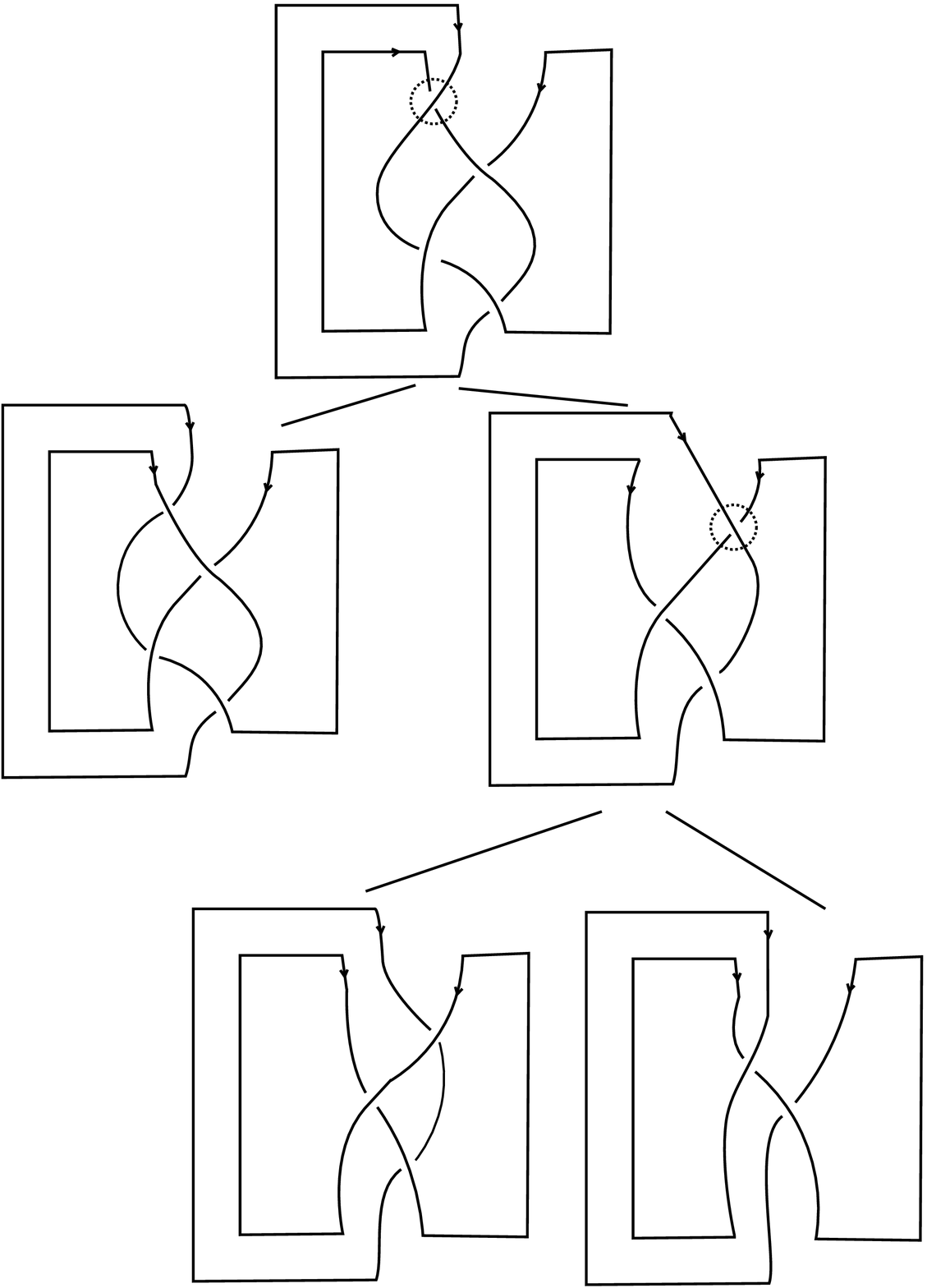,height=12cm}}
\centerline{\footnotesize{Fig. 1.4}}
As it is easily seen the leaves of the tree are trivial links and every branching reflects a certain operation on the diagram at the marked crossing point. 
To compute $w_L$ it is enough to have the following tree: \\ \ \\
\centerline{\psfig{figure=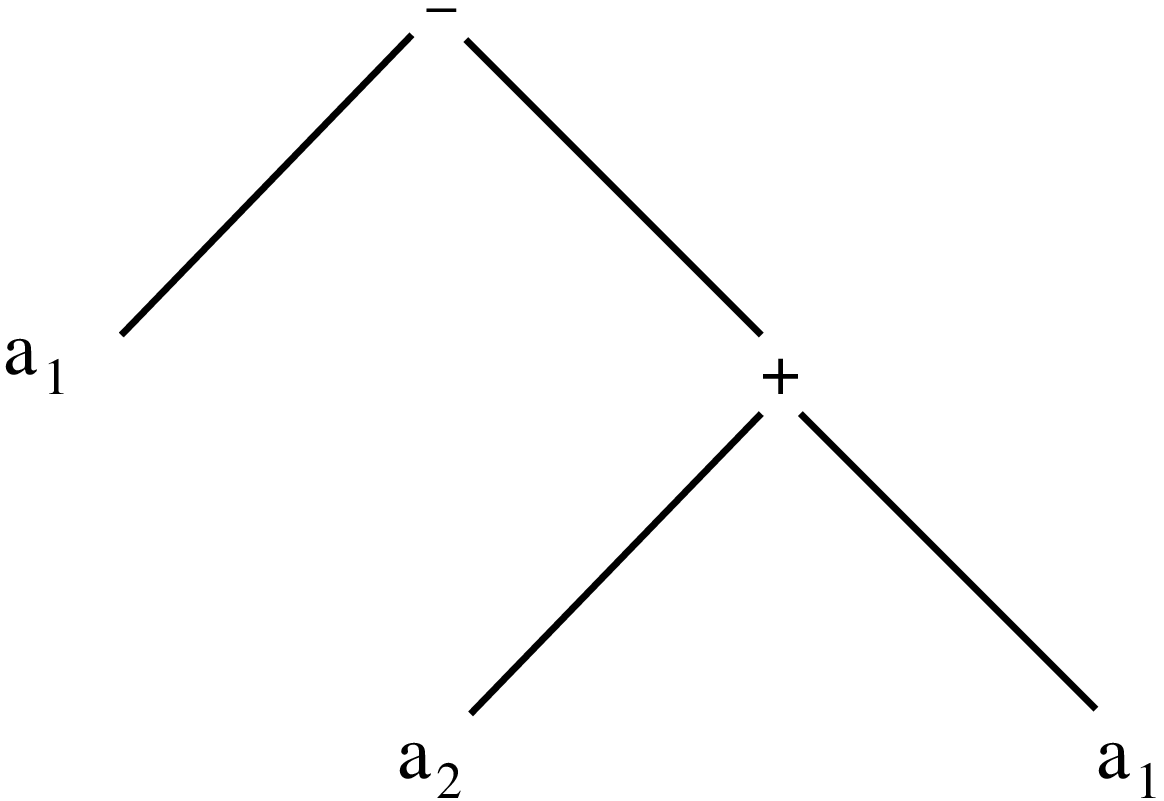,height=3cm}}
Here the sign indicates the sign of the crossing point at which the operation was performed, and the leaf entries are the 
values of $w$ for the resulting trivial links. Now we may conclude that
\[w_L\!=\!a_1\!*\!(a_2|a_1).\]
Such a binary tree of operations on the diagram resulting in trivial links at the leaves will be called the resolving 
tree of the diagram. \newline \indent There exists a standard procedure to obtain such a tree for every diagram. 
It will be described in the next paragraph and it will play an essential role in the proof of Theorem 1.8. 
It should be admitted that the idea is due to Ball and Metha {[$B$-$M$]} and we learned this from the Kauffman lecture notes {[$K$-3]}. \\

\section{Proof of the Main Theorem}\label{Section 2}
\begin{definition}\textbf{2.1.} Let $L$ be an oriented diagram of $n$ components and let $b\!=$($b,\dots,b_n$) be base points of $L$, one point 
from each component of $L$, but not the crossing points. Then we say that $L$ is untangled with respect to $b$ if the following holds: if one travels along $L$ 
(according to the orientation of $L$) starting from $b_1$, then, after having returned to $b_1-$ from $b_2,\dots,$ finally from $b_n$, 
then each crossing which is met for the first time is crossed by a bridge.\end{definition}
\indent It is easily seen that for every diagram $L$ of an oriented link there exists a resolving tree such that the leaf diagrams 
are untangled (with respect to appropriately chosen base points). This is obvious for diagrams with no crossings at all, and once it is known 
for diagrams with less than $n$ crossings we can use the following procedure for any diagram with $n$ crossings: choose base points arbitrarily 
and start walking along the diagram until the first ``bad" crossing $p$ is met, i.e. the first crossing which is crossed by a tunnel when first met. 
Then begin to construct the tree changing the diagram in this point. If, for example, sgn $p\!=\!+$ we get \\
\centerline{{\psfig{figure=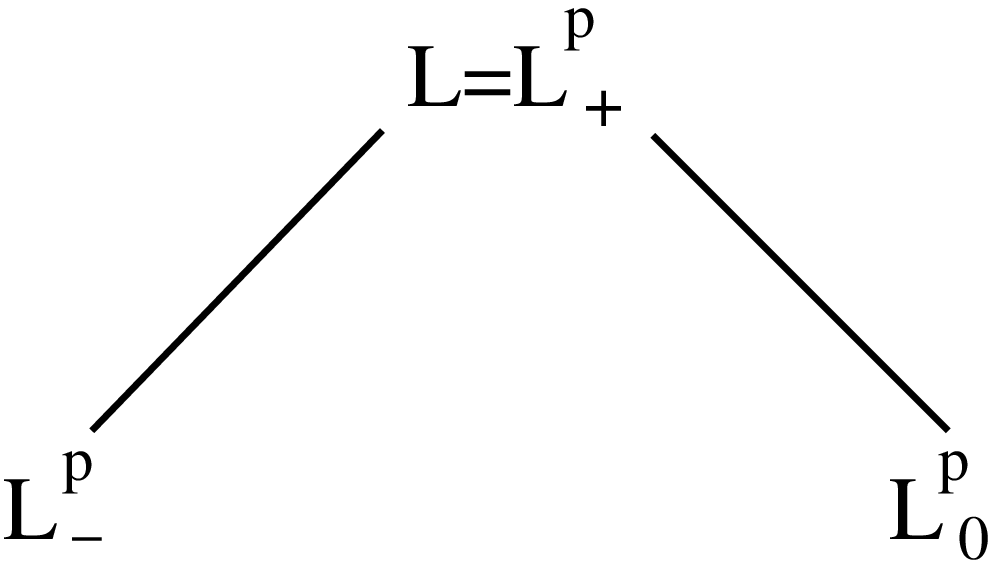,height=2.1cm}}}\ \\
Then we can apply the inductive hypothesis to $L_0^p$ and we can continue the procedure with $L_-^p$ (walking further along the diagram and looking 
for the next bad point). \\
\vspace{1mm}\\
\indent To prove Theorem 1.8 we will construct the function $w$ as defined on diagrams. In order to show that $w$ is an invariant 
of isotopy classes of oriented links we will verify that $w$ is preserved by the Reidemeister moves. \\
\vspace{1mm}\\
\indent We use induction on the number $cr(L)$ of crossing points in the diagram. For each $k \geqslant 0$ we define a function $w_k$ assigning an element 
of $\A$ to each diagram of an oriented link with no more than $k$ crossings. Then $w$ will be defined for every diagram by $w_L = w_k(L)$ 
where $k \geqslant cr(L)$. Of course the functions $w_k$ must satisfy certain coherence conditions for this to work. 
Finally we will obtain the required properties of $w$ from the properties of $w_k$'s. \\
\indent We begin from the definition of $w_0$. For a diagram $L$ of $n$ components with $cr(L)=0$ we put \\
\vspace{1mm}
\begin{flalign*}
    {2.2}&& w_0(L) &= a_n.&
\end{flalign*}
To define $w_{k+1}$ and prove its properties we will use the induction several times. To avoid misunderstandings the following will be called 
the ``Main Inductive Hypothesis": M.I.H. We assume that we have already defined a function $w_k$ attaching an element of $\A$ to each diagram $L$ 
for which $cr(L) \leqslant k$. We assume that $w_k$ has the following properties:
\begin{flalign*}
{2.3}&& w_k(U_n) &= a_n &
\end{flalign*}
for $U_n$ being an untangled diagram of $n$ components (with respect to some choice of base points).
\begin{flalign*}
{2.4} && w_k(L_+) &= w_k(L_-)|w_k(L_0)&\\
{2.5} && w_k(L_-) &= w_k(L_+)*w_k(L_0) &
\end{flalign*}
for $L_+$, $L_-$ and $L_0$ being related as usually.
\begin{flalign*}
{2.6}&& w_k(L) &= w_k(R(L))&
\end{flalign*}
where $R$ is a Reidemeister move on $L$ such that $cr(R(L))$ is still at most $k$. \\
\indent Then, as the reader may expect, we want to make the Main Inductive Step to obtain the existence of a function $w_{k+1}$ with analogous 
properties defined on diagrams with at most $k$ +1 crossings. \\
\indent Before dealing with the task of making the M.I.S. let us explain that it will really end the proof of the theorem. 
It is clear that the function $w_k$ satisfying M.I.H. is uniquely determined by properties 2.3, 2.4, 2.5 and the fact that for every diagram 
there exists a resolving tree with untangled leaf diagrams. Thus the compatibility of the function $w_k$ is obvious and they define 
a function $w$ defined on diagrams. \\
\indent The function $w$ satisfies conditions (2) and (3) of the theorem because the function $w_k$ satisfy such conditions. \\
\indent If $R$ is a Reidemeister move on a diagram $L$, then $cr(R(L))$ equals at most $k=cr(L)+2$, whence \\
\indent $w_{R(L)}$=$w_k(R(L))$, $w_L$=$w_k(L)$ and by the properties of $w_k, w_k(L)$=$w_k(R(L))$ what implies $w_{R(L))}$=$w_L$. 
It follows that $w$ is an invariant of the isotopy class of oriented links. \\
\indent Now it is clear that $w$ has the required property (1) too, since there is an untangled diagram $U_n$ in the same isotopy 
class as $T_n$ and we have $w_k(U_n)=a_n$. \\
\indent The rest of this section will be occupied by the M.I.S. For a given diagram $D$ with $cr(D)\!\leqslant\!k+1$ we will denote 
by $\mathcal{D}$ the set of all diagrams which are obtained from $D$ by operations of the kind 
\parbox{3.1cm}{\psfig{figure=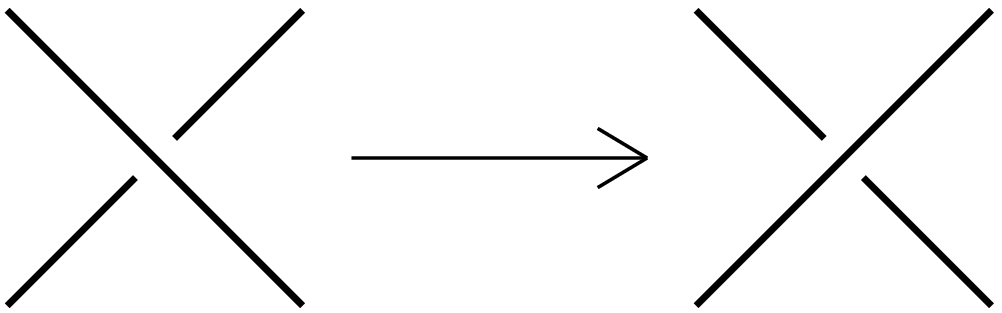,height=0.9cm}}
or
\parbox{3.1cm}{\psfig{figure=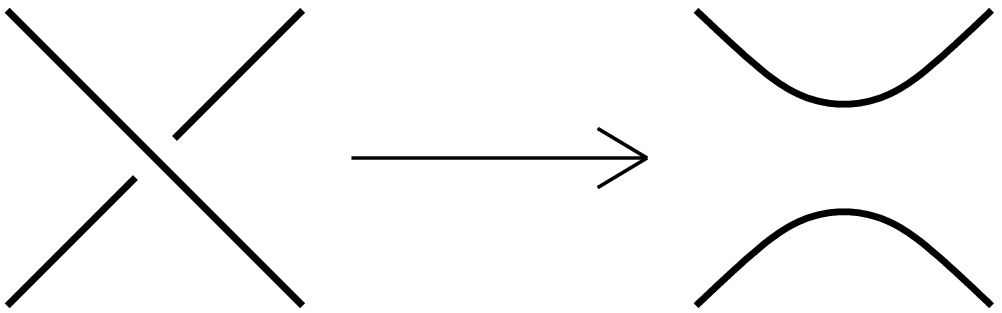,height=0.9cm}}. 
Of course, once base points $b$=($b_1,\dots, b_n$) are chosen on $D$, then the same points can be chosen as base points for 
any $L\!\in\!\mathcal{D}$, provided $L$ is obtained from $D$ by the operations of the first type only. \\
\indent Let us define a function $w_b$, for a given $D$ and $b$, assigning an element of $\A$ to each $L\!\in\!\mathcal{D}$. If $cr(L)\!<\!k+1$ we put
\begin{flalign*}
{2.7}&& w_b(L) &= w_k(L) &
\end{flalign*}
If $U_n$ is an untangled diagram with respect to $b$ we put
\begin{flalign*}
{2.8}&& w_b(U_n) &= a_{n} &
\end{flalign*}
($n$ denotes the number of components). \\
Now we can proceed by induction on the number $b(L)$ of bad crossings in $L$ (in the symbol $b(L)\ b$ works simultaneously 
for ``bad" and for $b$=($b_1,\dots,b_n$). For a different choice of base points $b'$=($b'_1,\dots,b'_n$) we will write $b'(L))$. 
Assume that $w_b$ is defined for all $L\!\in\!\mathcal{D}$ such that $b(L)\!<\!t$. Then for $L$, $b(L)$=$t$, let $p$ be the first 
bad crossing for $L$(starting from $b_1$ and walking along the diagram). Depending on $p$ being positive or negative we have $L$=$L_+^p$ or $L$=$L_-^p$. 
We put 
\begin{flalign*}
{2.9} & \ \ \ \ \ \ \ \ \ w_b(L)= 
\begin{cases}
    w_b(L_-^p) | w_b(L_0^p), & \text{if sgn } p = + \\
    w_b(L_+^p)*w_b(L_0^p), & \text{if sgn } p = -. \\
\end{cases}&
\end{flalign*}
We will show that $w_b$ is in fact independent of the choice of $b$ and that it has the properties required from $w_{k+1}$. \\
\vspace{1mm}\\
\textbf{Conway Relations for $\boldsymbol{w_b}$} \\
\vspace{1mm}\\
\indent Let us begin with the proof that $w_b$ has properties $2.4$ and $2.5$. We will denote by $p$ the considered crossing point. 
We restrict our attention to the case: $b(L_+^p)>b(L_-^p)$. The opposite situation is quite analogous. \\
\indent Now, we use induction on $b(L_-^p)$. If $b(L_-^p)$=$0$, then $b(L_+^p)$=$1$, $p$ is the only bad point of $L_+^p$, 
and by defining equalities $2.9$ we have
\[w_b(L_+^p)\!=\!w_b(L_-^p)|w_b(L_0^p)\]
and using 1.6 we obtain
\[w_b(L_-^p)\!=\!w_b(L_+^p)\!*\!w_b(L_0^p).\]
Assume now that the formulae $2.4$ and $2.5$ for $w_b$ are satisfied for every diagram $L$ such that $b(L_-^p)\!<\!t$, $t\!\geqslant\!1$. 
Let us consider the case $b(L_-^p)$=$t$. \\
\indent By the assumption $b(L_+^p)\!\geqslant\!2$. Let $q$ be the first bad point on $L_+^p$. Assume that $q$=$p$. Then by $2.9$ we have
\[w_b(L_+^p)\!=\!w_b(L_-^p)| w_b(L_0^p).\]
Assume that $q \neq p$. Let sgn $q=+$, for example. Then by 2.9 we have
\[w_b(L_+^p) = w_b(L{_+^p}{_+^q}) =   w_b(L{_+^p}{_-^q}) |  w_b(L{_+^p}{_ 0^q}).\]
But $b(L{_-^p}{_-^q})\!<\!t$ and $cr(L{_+^p}{_ 0^q})\!\leqslant\!k$, whence by the inductive hypothesis and M.I.H. we have
\[w_b(L{_+^p}{_-^q})\!=\!w_b(L{_-^p}{_-^q})|w_b(L{_0^p}{_ -^q})\]
and
\[w_b(L{_+^p}{_0^q})\!=\!w_b(L{_-^p}{_0^q})|w_b(L{_0^p}{_ 0^q})\]
whence
\[w_b(L_+^p) = (w_b(L{_-^p}{_-^q}) |  w_b(L{_0^p}{_-^q}))|( w_b(L{_-^p}{_ 0^q})| w_b(L{_0^p}{_0^q)})\]
and by the transposition property 1.3
\begin{flalign*}
{2.10} && w_b(L_+^p) &=(w_b(L{_-^p}{_-^q})|w_b(L{_-^p}{_0^q})) |(w_b(L{_0^p}{_ -^q})|w_b(L{_0^p}{_0^q})).&
\end{flalign*}
On the other hand $b(L{_-^p}{_-^q})\!<\!t$ and $cr(L_ 0^p)\!\leqslant\!k$, so using once more the inductive hypothesis and M.I.H. we obtain
\[
{2.11} \begin{split}
w_b(L_-^p) &= w_b(L{_-^p}{_+^q}) =   w_b(L{_-^p}{_-^q}) |  w_b(L{_-^p}{_ 0^q}) \\
w_b(L_0^p) &= w_b(L{_0^p}{_+^q}) =   w_b(L{_0^p}{_-^q}) |  w_b(L{_0^p}{_ 0^q})
\end{split}
\]
Putting 2.10 and 2.11 together we obtain
\[w_b(L_+^p) = w_b(L_-^p) |  w_b(L_0^p)\]
as required. If sgn $q=-$ we use $1.4$ instead of $1.3$. This completes the proof of Conway Relations for $w_b$. \\
\vspace{1mm}\\
\textbf{Changing Base Points} \\
\vspace{1mm}\\
\indent We will show now that $w_b$ does not depend on the choice of $b$, provided the order of components is not changed. It amounts to the 
verification that we may replace $b_i$ by $b'_i$ taken from the same component in such a way that $b'_i$ lies after $b_i$ and there is exactly 
one crossing point, say $p$, between $b_i$ and $b'_i$. Let $b'$=($b_1, \dots,b'_i,\dots,b_n$). We want to show that $w_b(L)\!=\!w_{b'}(L)$ for 
every diagram with $k+1$ crossings belonging to $\mathcal{D}$. We will only consider the case sgn $p=+$; the case sgn $p=-$ is quite analogous. \\
\indent We use induction on $B(L)=$max$(b(L),b'(L))$. We consider three cases. \\
\vspace{1mm}\\
\indent \textsc{Cbp 1.} Assume $B(L)\!=\!0$. Then $L$ is untangled with respect to both choices of base points and by 2.8
\[w_b(L) = a_n = w_{b'}(L).\]
\indent \textsc{Cbp 2.} Assume that $B(L)\!=\!1$ and $b(L)\!\neq\!b'(L)$. This is possible only when $p$ is a self-crossing point of the $i$-th component of $L$. There are two subcases to be considered. \\
\vspace{1mm}\\
\indent \textsc{Cbp 2} (a): $b(L)\!=\!1$ and $b'(L)\!=\!0$. Then $L$ is untangled with respect to $b'$ and by 2.8
$$  w_{b'}(L)=a_n $$
$$ \ \ \ \ w_b(L)=w_b(L_+^p)= w_b(L_-^p)|w_b(L_0^p)$$
\indent Again we have restricted our attention to the case sgn $p\!=\!+$.  Now, $w_b(L_-^p)\!=\!a_n$ since $b(L_-^p)\!=\!0$, and $L_0^p$ is 
untangled with respect to a proper choice of base points. Of course $L_0^p$ has $n+1$ components, so $w_b(L_0^p)\!=\!a_{n+1}$ by 2.8. 
It follows that $w_b(L)\!=\!a_n|a_{n+1}$ and $a_n|a_{n+1}=a_n$ by $1.1$. \\
\vspace{1mm}\\
\indent \textsc{Cbp 2}(b): $b(L)\!=\!0$ and $b'(L)\!=\!1$. This case can be dealt with like \textsc{Cbp 2}(a).\\
\vspace{1mm}\\
\indent \textsc{Cbp 3.} $B(L)\!=\!t\!>\!1$ or $B(L)\!=\!1\!=\!b(L)\!=\!b'(L)$. We assume by induction $w_b(K)\!=\!w_{b'}(K)$ for $B(K)\!<\!B(L)$. Let $q$ be a crossing point which is bad with respect to $b$ and $b'$ as well. We will consider this time the case sgn $q=-$. The case sgn $q=+$ is analogous. \\
\indent Using the already proven Conway relations for $w_b$ and $w_{b'}$ we obtain
\[w_b(L)= w_b(L_-^q)=w_b(L_+^q)\!*\!w_b(L_0^q) \]
\[w_{b'}(L)= w_{b'}(L_-^q)= w_{b'}(L_+^q)\!*\!w_{b'}(L_0^q)\]
But $B(L_+^q)\!<\!B(L)$ and $cr(L_0^q)\!\leqslant\!k$, whence by the inductive hypothesis and M.I.H. hold
\[w_b(L_+^q) = w_{b'}(L_+^q)\]
\[w_b(L_0^q) = w_{b'}(L_0^q)\]
which imply $w_b(L)\!=\!w_{b'}(L)$. This completes the proof of this step (C.B.P.). \\
\indent Since $w_b$ turned out to be independent of base point changes which preserve the order of components we can now consider a function $w^0$ to be defined in such a way that it attaches an element of $\A$ to every diagram $L$, $cr(L) \leqslant k+1$ with a fixed ordering of components. \\
\vspace{1mm}\\
\textbf{Independence of $\boldsymbol{w^0}$ of Reidemeister Moves} (I.R.M) \\
\vspace{1mm}\\
\indent When $L$ is a diagram with fixed order of components and $R$ is a Reidemeister move on $L$, then we have a natural ordering of components on $R(L)$. 
We will show now that $w^0(L)=w^0(R(L))$. Of course we assume that $cr(L)$, $cr(R(L)) \leqslant k+1$. \\
\indent We use induction on $b(L)$ with respect to properly chosen base points $b=(b_1, \dots , b_n)$. 
Of course the choice must be compatible with the given ordering of components. 
We choose the base points to lie outside the part of the diagram involved in the considered Reidemeister move $R$, 
so that the same points may work for the diagram $R(L)$ as well. We have to consider the three standard types of Reidemeister moves (Fig. 2.1). \\ \ \\
\centerline{\psfig{figure=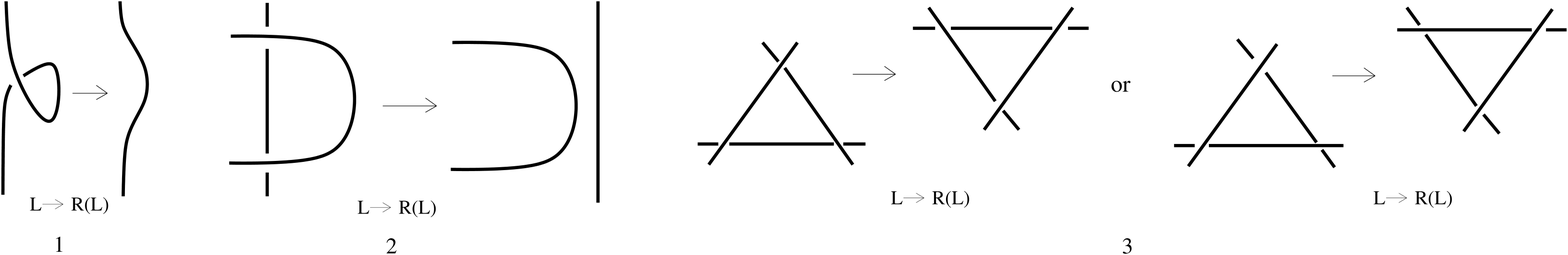,height=2.5cm}}
\centerline{\footnotesize{Fig. 2.1}}
\indent Assume that $b(L)=0$. Then it is easily seen that also $b(R(L))=0$, and the number of components is not changed. Thus
\[w^0(L)=w^0(R(L))\ \textup{ by\ 2.8.}\]
\indent We assume now by induction that $w^0(L)=w^0(R(L))$ for $b(L)<t$. Let us consider the case $b(L)=t$. Assume that there is a bad crossing $p$ in $L$ which is different from all the crossings involved in the considered Reidemeister move. Assume, for example, that sgn $p=+$. Then, by the inductive hypothesis, we have
\begin{flalign*}
    {2.12} && w^0(L_-^p) &= w^0(R(L_-^p))&
\end{flalign*}
and by M.I.H.
\begin{flalign*}
    {2.13} && w^0(L_0^p) &= w^0(R(L_0^p))&
\end{flalign*}
Now, by the Conway relation 2.4, which was already verified for $w^0$ we have
\[w^0(L)=w^0(L_+^p)=w^0(L_-^p)|w^0(L_0^p) \]
\[w^0(R(L))=w^0(R(L)_+^p)=w^0(R(L)_-^p)|w^0(R(L)_0^p) \]
whence by 2.12 and 2.13
\[w^0(L)=w^0(R(L))\]
Obviously $R(L_-^p)=R(L)_-^p$ and $R(L_0^p)=R(L)_0^p$. \\
\indent It remains to consider the case when $L$ has no bad points, except those involved in the considered Reidemeister move. 
We will consider the three types of moves separately. 
The most complicated is the case of a Reidemeister move of the third type. To deal with it let us formulate the following observation: \\
\indent Whatever the choice of base points is, the crossing point of the top arc and the bottom arc cannot be the only bad point of the diagram. \\
\centerline{\psfig{figure=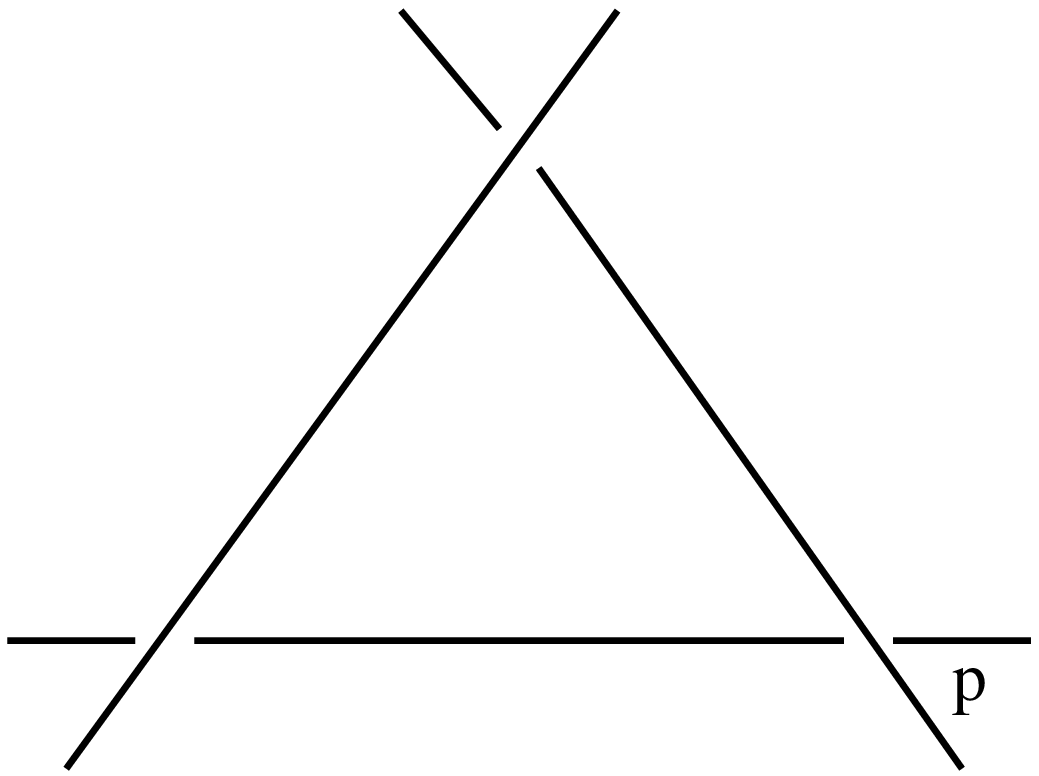,height=4.1cm}}
\centerline{\footnotesize{Fig. 2.2}}
The proof of the above observation amounts to an easy case by case checking and we omit it. The observation makes possible the following 
induction: we can assume that we have a bad point at the crossing between the middle arc and the lower or the upper arc. 
Let us consider for example the first possibility; thus $p$ from Fig. 2.2 is assumed to be a bad point. 
We consider two subcases, according to sgn $p$ being $+$ or $-$. \\
\indent Assume sgn $p=+$. Then by Conway relations
\[w^0(L)=w^0(L_+^p)=w^0(L_-^p)|w^0(L_0^p)\]
\[w^0(R(L))=w^0(R(L)_+^p)=w^0(R(L)_-^p)|w^0(R(L)_0^p)\]
But $R(L)_-^p=R(L_-^p)$ and by the inductive hypothesis
\[w^0(L_-^p)=w^0(R(L_-^p))\]
Also $R(L)_0^p$ is obtained from $L_0^p$ by two subsequent Reidemeister moves of type two (see Fig. 2.3), whence by M.I.H.
\[w^0(R(L)_0^p)=w^0(L_0^p)\]
and the equality
\[w^0(L)=w^0(R(L))\ \textup{follows.}\]

\centerline{\epsfig{figure=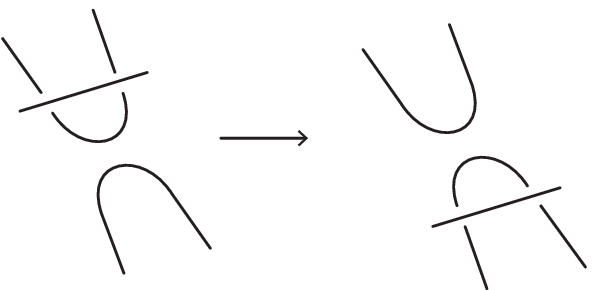,height=3.5cm}}
\centerline{\footnotesize{Fig.2.3}}
Assume now that sgn $p=-$. Then by Conway relations
\[w^0(L)=w^0(L_-^p)=w^0(L_+^p)*w^0(L_0^p) \]
\[w^0(R(L))=w^0(R(L)_-^p)=w^0(R(L)_+^p)*w^0(R(L)_0^p)\]
But $R(L)_+^p=R(L_+^p)$ and by the inductive hypothesis
\[w^0(L_+^p)=w^0(R(L_+^p))\]
Now, $L_0^p$ and $R(L)_0^p$ are essentially the same diagrams (see Fig. 2.4), whence $w^0(L_0^p)=w^0(R(L)_0^p)$ and the equality
\[w^0(L)=w^0(R(L))\ \textup{follows.}\]
\centerline{\epsfig{figure=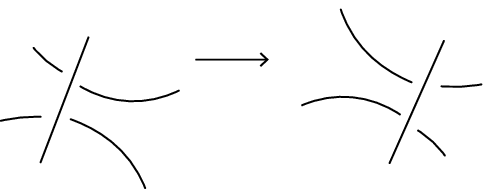, height=3.5cm}}
\centerline{\footnotesize{Fig. 2.4}} 
\ \\ \ \\
Reidemeister moves of the first type. The base points can always be chosen so that the crossing point involved in the move is good. \\
Reidemeister moves of the second type. There is only one case, when we cannot choose base points to guarantee the points involved in the move to be good. 
It happens when the involved arcs are parts of different components and the lower arc is a part of the earlier component. 
In this case the both crossing points involved are bad and they are of different signs, of course. Let us consider the situation shown in Fig. 2.5. \\
\centerline{\epsfig{figure=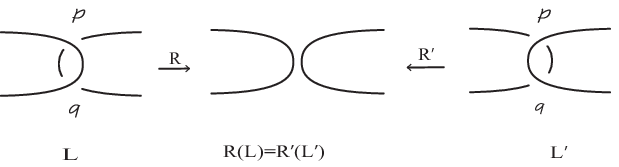, height=3.5cm}}
\centerline{\footnotesize{Fig. 2.5}}
We want to show that $w^0(R(L))=w^0(L)$. But by the inductive hypothesis we have
\[w^0(L')=w^0(R'(L'))=w^0(R(L)).\]
Using the already proven Conway relations, formulae 1.6 and 1.7 and M.I.H. if necessary, it can be proved that $w^0(L)=w^0(L')$. 
Let us discuss in detail the case involving M.I.H. It occurs when sgn $p=-$. 
Then we have \[w^0(L)=w^0(L_+^q)=w^0(L_-^q)|w^0(L_0^q)=(w^0(L{_-^q}{_+^p})*w^0(L{_-^q}{_0^p}))|w^0(L_0^q)\]
But $L{_-^q}{_+^p}= L'$ and by M.I.H. $w^0(L{_-^q}{_0^p})=w^0(L_0^q)$ (see Fig. 2.6, 
where $L{_-^q}{_0^p}$ and $L_0^q$ are both obtained from $K$ by a Reidemeister move of the first type). \\ \ \\
\centerline{\epsfig{figure=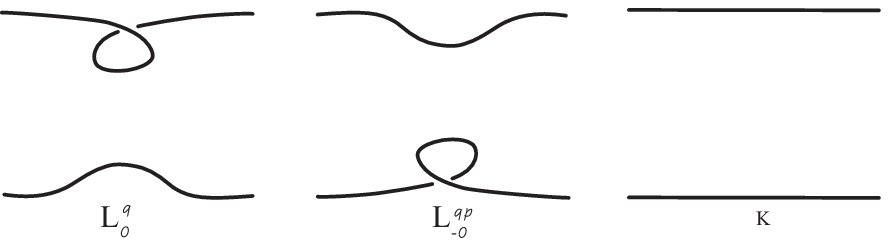, height=3.2cm}}
\centerline{\footnotesize{Fig. 2.6}}
Thus by 1.7:
\[w^0(L) = w^0(L')\ \textup{whence}\]
\[w^0(L)=w^0(R(L)). \]
The case: sgn $p=+$ is even simpler and we omit it. This completes the proof of the independence of $w^0$ of Reidemeister moves. 
To complete the Main Inductive Step it is enough to prove the independence of $w^0$ of the order of components. 
Then we set $w_k = w^0$. The required properties have been already checked. \\
\vspace{1mm}\\
\textbf {Independence of the Order of Components} (I.O.C.) \\
\vspace{1mm}\\
\indent It is enough to verify that for a given diagram $L\ (cr(L) \leqslant k+1)$ and fixed base points $b=(b_1, \dots , b_i, b_{i+1}, \dots , b_n)$ 
we have \[w_b(L)=w_{b'}(L)\]
where $b'=(b_1, \dots , b_{i+1}, b_{i}, \dots , b_n)$. This is easily reduced by the usual induction on $b(L)$ to the case of an untangled diagram. 
To deal with this case we will choose $b$ in an appropriate way. \\
\indent Before we do it, let us formulate the following observation: If $L_i$ is a trivial component of $L$, i.e. $L_i$ has no crossing points, 
neither with itself, nor with other components, then the specific position of $L_i$ in the plane has no effect on $w^0(L)$; 
in particular we may assume that $L_i$ lies separately from the rest of the diagram: \\ \ \\
\centerline{\epsfig{figure=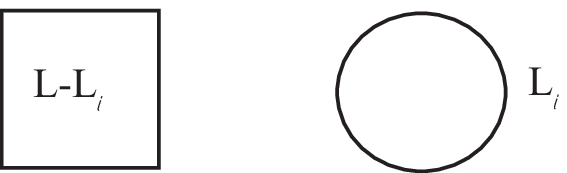, height=2.2cm}}
\centerline{\footnotesize{Fig. 2.7}}
This can be easily achieved by induction on $b(L)$, or better by saying that it is obvious. \\
\indent For an untangled diagram we will be done if we show that it can be transformed into another one with less crossings by a series of 
Reidemeister moves which do not increase the crossing number. We can then use I.R.M. and M.I.H. This is guaranteed by the following lemma.
\begin{lemma}\textbf{2.14.} Let $L$ be a diagram with $k$ crossings and a given ordering of components $L_1, L_2, \dots , L_n$. 
Then either $L$ has a trivial circle as a component or there is a choice of base points $b=(b_1, \dots , b_n)$; $b_i \in L_i$ such 
that an untangled diagram $L^{u}$ associated with $L$ and $b$ (that is all the bad crossings of $L$ are changed to good ones) 
can be changed into a diagram with less than $k$ crossings by a sequence of Reidemeister moves not increasing the number of crossings.
\end{lemma}
This was probably known to Reidemeister already, however we prove it in the Appendix for the sake of completeness.\\
\indent With I.O.C. proven we have completed M.I.S. and the proof of Theorem 1.8.\\
\section{Quasi-algebras}\label{Section 3}
\indent We shall now describe a certain generalization of Theorem 1.8. This is based on the observation that it was not necessary 
to have the operations $|$ and $*$ defined on the whole product $\A \times \A$. Let us begin with the following definition. 
\begin{definition}\textbf{3.0.}
A quasi Conway algebra is a triple $(\A, B_1, B_2), B_1, B_2$ being subsets of $\A \times \A$, together with 0-argument 
operations $a_1,a_2, \dots ,a_n, \dots$ and two 2-argument operations $|$ and $*$ defined on $B_1$ and $B_2$ respectively satisfying the conditions:
\end{definition}
\begin{align*}
\left. \begin{aligned}
3.1 && \qquad &a_n | a_{n+1} = a_n \\
3.2 && \qquad &a_n * a_{n+1} = a_n \\
3.3&& \qquad &(a|b)|(c|d) = (a|c)|(b|d) \\
3.4 && \qquad &(a|b)*(c|d) = (a*c)|(b*d) \\
3.5 && \qquad &(a*b)*(c*d) = (a*c)*(b*d) \\
3.6 && \qquad &(a|b)*b = a \\
3.7 && \qquad &(a*b)|b = a. \\
\end{aligned}
\right\}
\indent \text{whenever the both
sides are defined}
\end{align*}
\indent We would like to construct invariants of Conway type using such quasi-algebras. As before $a_n$ will be the value of 
the invariant for the trivial link of $n$ components. \\
\indent We say that $\A$ is geometrically sufficient if and only if for every resolving tree of each diagram of an oriented link all the operations 
that are necessary to compute the root value are defined.
\begin{theorem}\textbf{3.8.} Let $\A$ be a geometrically sufficient quasi Conway algebra. 
There exists a unique invariant $w$ attaching to each isotopy class of links an element of $\A$ and satisfying the conditions
\end{theorem}
\begin{enumerate}
\item $w_{T_n}=a_n$ for $T_n$ being a trivial link of n components,\\
\item if $L_+$, $L_-$ and $L_0$ are diagrams from Fig. 1.1, then
\[w_{L_+}=w_{L_-}|w_{L_0}  \text{ and } \]
\[w_{L_-}=w_{L_+}*w_{L_0}.\]
\end{enumerate}
\noindent The proof is identical with the proof of Theorem 1.8.\\
\indent As an example we will now describe an invariant, whose values are polynomials in an infinite number of variables.
\begin{example}\textbf{3.9.} $\A=N \times Z[ y_1^{\mp1},{x'}_2^{\mp 1},{z'}_2,\ x_1^{\mp1},z_1,\ x_2^{\mp1},z_2,x_3^{\mp1},z_3,\dots], 
B_1=B_2=B=\{ ((n_1,\ w_1),(n_2,w_2))\in\A\times\A:|n_1-n_2|=1\}, a_1=(1,1),a_2=(2, x_1+y_1+z_1),\dots,
a_n=(n,\Pi_{i=1}^{n-1}(x_i+y_i)+z_1\Pi_{i=2}^{n-1} (x_i+y_i)+\dots+z_{n-2}(x_{n-1}+y_{n-1})+z_{n-1}),\dots$ 
where $y_i=x_i\frac{y_1}{x_1}$. To define the operations $|$ and $*$ consider the following system of equations:\end{example}
\begin{flalign*}
{(1)} &&  &x_1w_1+y_1w_2=w_0-z_1& \\
{(2)} &&  &x_2w_1+y_2w_2=w_0-z_2& \\
{(2^\prime)}&& &x'_2w_1+y'_2w_2=w_0-z'_2& \\
{(3)} && &x_3w_1+y_3w_2=w_0 -z_3 \\
{(3^\prime)}&& &x'_3w_1+y'_3w_2=w_0-z'_3 \\
&& &\dots \\
{(i)}&& &x_iw_1+y_iw_2=w_0-z_i \\
{\text({i}^\prime})&& &x'_iw_1+y'_iw_2=w_0-z'_i \\
&& &\ldots
\end{flalign*}
where $y'_i=\frac{x'_iy_1}{x_i},x'_i=\frac{x'_2x_1}{x_{i-1}}$ and $z'_i$ are defined inductively to satisfy
\[\frac{z'_{i+1}-z_{i-1}}{x_1x'_2}= \Big(1+\frac{y_1}{x_1}\Big)\Big(\frac{z'_i}{x'_i}-\frac{z_i}{x_i}\Big).\]
We define $(n,w)=(n_1,w_1)|(n_2,w_2)$ (resp.$(n,w)=(n_1,w_1)*(n_2, w_2))$ as follows: $n=n_1$ and if $n_1=n_2-1$ 
then we use the equations $(n)$ to get $w$; namely $x_nw+y_nw_1=w_2-z_n$ (resp. $x_nw_1+y_nw=w_2-z_n$). 
If $n_1=n_2+1$ then we use the equation ($n'$) to get $w$; namely $x'_nw+y'_nw_1=w_2-z'_n$ (resp. $x'_nw_1+y'_nw=w_2-z'_n$). 
We can think of Example 1.11 as being a special case of Example 3.9. \\
\indent Now we will show that the quasi-algebra $\A$ for Example 3.9 satisfies the relations 1.1-1.7.\\
\indent It is an easy task to check that the first coordinate of elements from $\A$ satisfies the relations 1.1-1.7 
(compare with Example 1.9) and to check the relations 1.1, 1.2, 1.6 and 1.7 so we will concentrate our attention on relations 1.3, 1.4, and 1.5.\\
\indent It is convenient to use the following notation: if $w\in\A$ then $w=(\lvert{w}\rvert,F)$ and for
\[w_1|w_2=(\lvert{w_1}\rvert,F_1)|(\lvert{w_2}\rvert, F_2)=(\lvert{w}\rvert,F)=w \]
to use the notation
\begin{equation*} F=\begin{cases}
    F_1|_nF_2 \ \textup{if} \ n=\lvert{w_1}\rvert=\lvert{w_2}\rvert-1 \\
    F_1|_{n'}F_2 \ \textup{if} \ n=\lvert{w_1}\rvert=\lvert{w_2}\rvert+1.
\end{cases}
\end{equation*}
Similar notation we use for the operation $*$. \\

\indent In order to verify relations 1.3-1.5 we have to consider three main cases:\\
$1. \quad \lvert{a}\rvert=\lvert{c}\rvert-1=\lvert{b}\rvert+1=n$ \\
Relations 1.3-1.5 make sense iff $\lvert{d}\rvert=n$. The relation 1.3 has the form:
\[(F_a|_{n'}F_b)|_n(F_c|_{(n+1)'}F_d)=(F_a|_{n}F_c)|_{n'}(F_b|_{(n-1)}F_d).\]
From this we get:
\[\frac{1}{x_nx'_{n+1}}F_d-\frac{y'_{n+1}}{x_nx'_{n+1}}F_c-\frac{y_n}{x_nx'_n}F_b + \frac{y_ny'_n}{x_nx'_n}
F_a-\frac{z'_{n+1}}{x_nx'_{n+1}}-\frac{z_n}{x_n}+ 
\frac{y_nz'_n}{x_nx'_n}= \]
\[= \frac{1}{x'_nx_{n-1}} F_d-\frac{y_{n-1}}{x'_nx_{n-1}}F_b-\frac{y'_n}{x_nx'_n}F_c + \frac{y_ny'_n}{x_nx'_n}F_a-\frac{z_{n-1}}{x'_nx_{n-1}}-\frac{z'_n}{x'_n}+
\frac{y'_nz_n}{x_nx'_n}\]
Therefore:
\begin{align*}
\text{(i)}&& &x_{n-1}x'_n=x_nx'_{n+1}\\
\text{(ii)}&& &\frac{y'_{n+1}}{x'_{n+1}}=\frac{y'_n}{x'_n} \\
\text{(iii)}&& &\frac{y_n}{x_n}=\frac{y_{n-1}}{x_{n-1}} \\
\text{(iv)}&& &\frac{z'_{n+1}}{x_nx'_{n+1}}+\frac{z_n}{x_n}-\frac{y_nz'_n}{x_nx'_n}=\frac{z_{n-1}}{x'_nx_{n-1}}+\frac{z'_n}{x'_n}-\frac{y'_nz_n}{x_nx'_n}
\end{align*}
When checking the relations 1.4 and 1.5 we get exactly the same conditions (i)-(iv).\\
$2. \quad \lvert a \rvert =\lvert b \rvert -1 = \lvert{c} \rvert -1 = n$.\\
\indent (I) $\lvert{d}\rvert=n.$\\
The relation 1.3 has the following form:
\[(F_a|_{n}F_b)|_n(F_c|_{(n+1)'}F_d)=(F_a|_{n}F_c)|_{n}(F_b|_{(n+1)'}F_d).\]
We get after some calculations that it is equivalent to
\begin{align*}
\text{(v)}&& &\frac{y_n}{x_n}=\frac{y'_{n+1}}{x'_{n+1}}&
\end{align*}
The relations 1.4 and 1.5 reduce to the same condition (v). \\
\indent (II) $\lvert{d}\rvert=n+2.$ \\
Then the relations 1.3-1.5 reduce to the condition (iii). \\
$3. \quad \lvert{a}\rvert=\lvert{b}\rvert+1=\lvert{c}\rvert+1=n$ \\
\indent (I) $\lvert{d}\rvert = n-2$ \\
\indent (II) $\lvert{d}\rvert = n.$ \\
We get, after some computations, that relations 3 (I) and 3 (II) follow from the conditions (iii) and (v). \\
\indent Conditions (i)-(v) are equivalent to the conditions on $x'_i, y_i, y'_i$ and $z'_i$ described in Example 3.9. 
Therefore the quasi-algebra $\A$ from Example 3.9 satisfies the relations 1.1-1.7. Furthermore, if $L$ is a diagram and $p-$ its crossing, 
then the number of components of $L_0^p$ is always equal to the number of components of $L$ plus or minus one, so the set $B\subset\A\times\A$ 
is sufficient to define the link invariant associated with $\A$.\\

\section{Final remarks and problems}\label{Section 4}
\begin{remark}\textbf{4.1.} Each invariant of links can be used to build a better invariant which will be called weighted simplex of the invariant. 
Namely, if $w$ is an invariant and $L$ is a link of $n$ components $L_1,\dots, L_n$ then we consider an $n-1$ dimensional simplex 
$\Delta^{n-1}=(q_1,\dots,q_n)$. We associate with each face ($q_{i_1},\dots,q_{i_k}$) of $\Delta^{n-1}$ the value $w_{L'}$ where 
$L' = L_{i_1}\cup \cdots \cup L_{i_k}$.
\end{remark}
\indent We say that two weighted simplices are equivalent if there exists a bijection of their vertices which preserves weights of faces. 
Of course, the weighted simplex of an invariant of isotopy classes of oriented links is also an invariant of isotopy classes of oriented links. \\
\indent Before we present some examples, we will introduce an equivalence relation $\thicksim_c$ (Conway equivalence relation) on 
isotopy classes of oriented links ($\mathcal{L}$).
\begin{definition}\textbf{4.2.} $\thicksim_c$ is the smallest equivalence relation on $\mathcal{L}$ which satisfies the following condition: 
let $L'_1$ (resp. $L'_2$) be a diagram of a link $L_1$ (resp. $L_2$) with a given crossing $p_1$ (resp. $p_2$) such that $p_1$ and $p_2$ 
are crossings of the same sign and 
\[
\begin{split}
    (L'_1)_-^{p_1} \thicksim_c (L'_2)_-^{p_2}\ &\text{and} \\
    (L'_1)_0^{p_1} \thicksim_c (L'_2)_0^{p_2}
\end{split}
\]
then $L_1\thicksim_c L_2$.\\
\end{definition}
\indent It is obvious that an invariant given by a quasi Conway algebra is a Conway equivalence invariant.
\begin{example}\textbf{4.3.}{(a)} Two links shown on Fig. 4.1 are Conway equivalent but they can be distinguished by weighted simplices 
of the linking numbers. \end{example}
\centerline{\epsfig{figure= 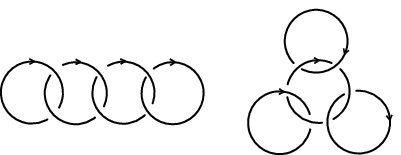,height=3.5cm}}
\centerline{\footnotesize{Fig. 4.1}}
{(b)} J. Birman has found three-braids (we use the notation of [M]);
\[
\begin{split}
    \gamma_1=\sigma_1^{-2}\sigma_2^3\sigma_1^{-1}\sigma_2^4\sigma_1^{-2}\sigma_2^{4}\sigma_1^{-1}\sigma_2 \\
    \gamma_2=\sigma_1^{-2}\sigma_2^3\sigma_1^{-1}\sigma_2^4\sigma_1^{-1}\sigma_2\sigma_1^{-2}\sigma_2^4
\end{split}
\]
which closures have the same values of all invariants described in our examples and the same signature but which can be distinguished by 
weighted simplices of the linking numbers [B]. \\
\indent As the referee has kindly pointed out the polynomial invariants described in 1.11 (a) and 1.11 (b) are equivalent. 
Namely if we denote them by $w_L$ and $w'_L$ respectively, then we have
\[w_L(x,\ y,\ z)=\Big(1-\frac{z}{1-x-y}\Big)w'_L(x,\ y)+\frac{z}{1-x-y}.\]
\begin{problem}\textbf{4.4.}
\begin{enumerate}
    \item[(a)]
 Is the invariant described in Example 3.9 better than the polynomial invariant 
from Example 1.11?\footnote{Added for e-print: Adam Sikora proved in his
Warsaw master degree thesis written under direction of P.Traczyk, that
the answer to Problem 4.4 (a) is negative, \cite{Si-1}.}
    \item[(b)]  Find an example of links which have the same polynomial invariant of Example 3.9 but which can be distinguished by some invariant 
given by a Conway algebra.\footnote{Added for e-print: Adam Sikora proved no invariant coming from a Conway algebra can distinguish links with the same
 polynomial invariant of Example 1.11, \cite{Si-2}.}
    \item[(c)] Do there exist two links $L_1$ and $L_2$ which are not Conway equivalent but which cannot be distinguished using any Conway algebra?
    \item[(d)]  Birman [$B$] described two closed 3-braid knots given by
\[
\begin{split}
    y_1=\sigma_1^{-3}\sigma_2^4\sigma_1^{-1}\sigma_2^5\sigma_1^{-3}\sigma_2^{5}\sigma_1^{-2}\sigma_2 \\
    y_2=\sigma_1^{-3}\sigma_2^4\sigma_1^{-1}\sigma_2^5\sigma_1^{-2}\sigma_2\sigma_1^{-3}\sigma_2^5
\end{split}
\]
which are not distinguished by the invariants described in our examples and by the signature. Are they Conway equivalent? 
(they are not isotopic because their incompressible surfaces are different).
\end{enumerate}
\end{problem}

\begin{problem}\textbf{4.5.} Given two Conway equivalent links, do they necessarily have the same signature?
\end{problem}
\indent The examples of Birman and Lozano [$B$; Prop. 1 and 2] have different signature but the same polynomial invariant of Example 1.11 (b)(see [$B$]).
\begin{problem}\textbf{4.6.} Let ($V_1$, $V_2$) be a Heegaard splitting of a closed 3-manifold $M$. Is it possible to modify the above 
approach using projections of links onto the Heegaard surface $\p V_1$?
\end{problem}
\indent We have obtained the results of this paper in early December 1984 and we were not aware at the time that an important part of our results 
(the invariant described in Example 1.11 (b)) had been got three months before us by four groups of researchers: 
R. Lickorish and K. Millett, J. Hoste, A. Ocneanu, P. Freyd and D. Yetter and that the first two groups used arguments similar to ours. 
We have been informed about this by J. Birman (letter received on January 28, '85) and by J. Montesinos (letter received on February 11, '85; 
it included the paper by R. Lickorish and K. Millett and also a small joint paper by all above mentioned mathematicians). \\

\section{Appendix}
Here we prove Lemma 2.14.\\
\indent A closed part cut out of the plane by arcs of $L$ is called an $i$-gon if it has $i$ vertices (see Fig. 5.1). \\
\centerline{{\psfig{figure=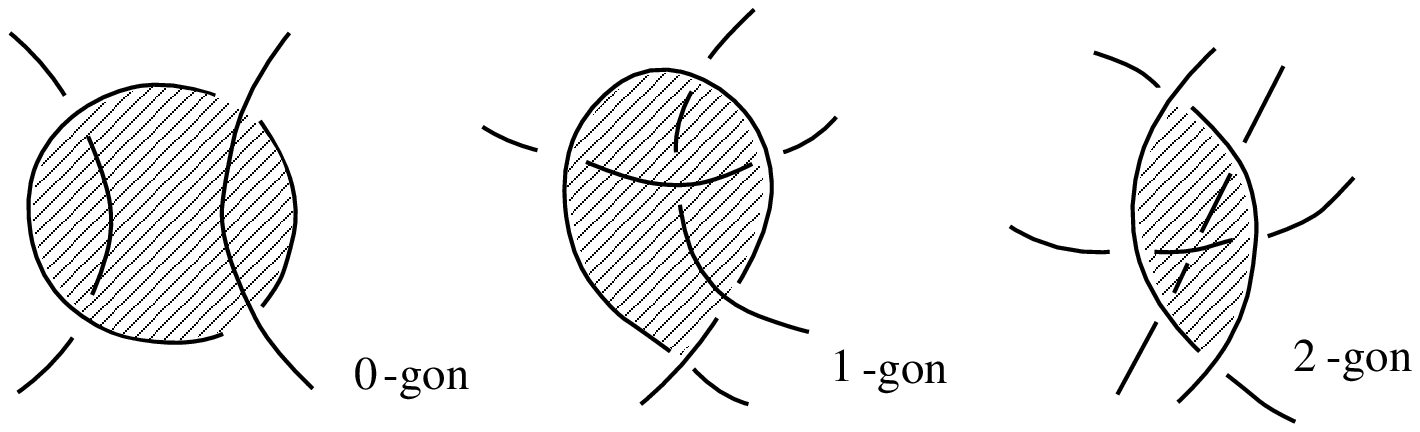,height=4.0cm}}}\ \\
\centerline{\footnotesize{Fig. 5.1}}
Every $i$-gon with $i\leqslant2$ will be called $f$-gon ($f$ works for few). Now let $X$ be an innermost $f$-gon that is an $f$-gon which 
does not contain any other $f$-gon inside. \\
\indent If $X$ is 0-gon we are done because $\p X$ is a trivial circle. If $X$ is 1-gon then we are done because int $X \cap L = \emptyset$ 
so we can perform on $L^u$ a Reidemeister move which decreases the number of crossings of $L^u$ (Fig. 5.2). \\
\centerline{{\psfig{figure=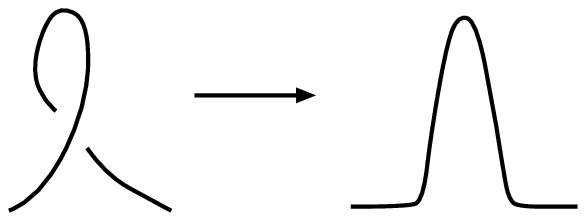,height=3.5cm}}}\ \\
\centerline{\footnotesize{Fig. 5.2}}\\
Therefore we assume that $X$ is a 2-gon. Each arc which cuts int $X$ goes from one edge to another. Furthermore, no component of $L$ lies 
fully in $X$ so we can choose base points $b=(b_1,\dots,b_n)$ lying outside $X$. This has important consequences: if $L^u$ is an untangled 
diagram associated with $L$ and $b$ then each 3-gon in $X$ supports a Reidemeister move of the third type (i.e. the situation of the Fig. 5.3 is impossible).\\
\centerline{{\psfig{figure=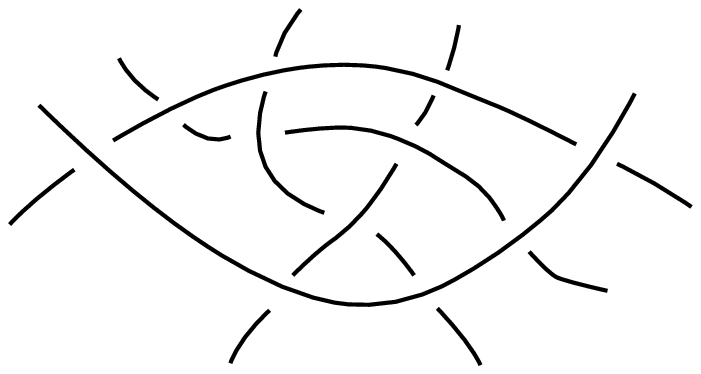,height=4.5cm}}}\ \\
\centerline{\footnotesize{Fig. 5.3}}
\indent Now we will prove Lemma 2.14 by induction on the number of crossings of $L$ contained in the 2-gon $X$ (we denote this number by $c$). \\
\indent If $c=2$ then int $X \cap L=\emptyset$ and we are done by the previous remark (2-gon $X$ can be used to make the Reidemeister move 
of the second type on $L^u$ and to reduce the number of crossings in $L^u$). \\
\indent Assume that $L$ has $c>2$ crossings in $X$ and that Lemma 2.14 is proved for less than $c$ crossings in $X$. 
In order to make the inductive step we need the following fact. 
\begin{proposition}\textbf{5.1.} If $X$ is an innermost 2-gon with int $X\cap L \neq \emptyset$ then there is a 3-gon $\Delta \subset X$ 
such that $\Delta \cap \p X \neq \emptyset$, int $\Delta \cap L = \emptyset$.
\end{proposition}
Before we prove Proposition 5.1 we will show how Lemma 2.14 follows from it. \\
\indent We can perform the Reidemeister move of the third type using the 3-gon $\Delta$ and reduce the number of 
crossings of $L^u$ in $X$ (compare Fig. 5.4).\\
\centerline{\epsfig{figure= 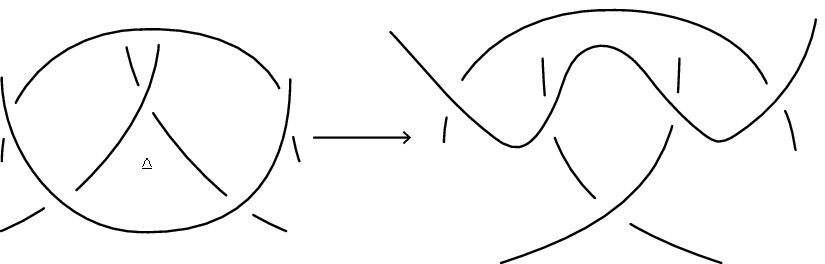, height=3.5cm}} 
\centerline{\footnotesize{Fig. 5.4}}\\
Now either $X$ is an innermost $f$-gon with less than $c$ crossings in $X$ or it contains an innermost $f$-gon with less that $c$ crossings in it. 
In both cases we can use the inductive hypothesis. \\ \indent Instead of proving Proposition 5.1 we will show a more general fact, 
which has Proposition 5.1 as a special case.
\begin{proposition}\textbf{5.2.} Consider a 3-gon $Y=$($a$, $b$, $c$) such that each arc which cuts it goes from the edge $\overline{ab}$ to the 
edge $\overline{ac}$ without self-intersections (we allow $Y$ to be a 2-gon considered as a degenerated 3-gon with $\overline{bc}$ collapsed 
to a point). Furthermore let int $Y$ be cut by some arc. Then there is a 3-gon $\Delta \subset Y$ such that $\Delta \cap \overline{ab} \neq\emptyset$ 
and int $\Delta$ is not cut by any arc.
\end{proposition}
\indent \textsc{Proof of Proposition 5.2:} We proceed by induction on the number of arcs in int $Y\cap L$ (each such an arc cuts $\overline{ab}$ 
and $\overline{ac}$). For one arc it is obvious (Fig. 5.5). Assume it is true for $k$ arcs ($k\geqslant 1$) and consider ($k+1)$-th arc $\gamma$. 
Let $\Delta_0=$($a_1$, $b_1$, $c_1$) be a $3$-gon from the inductive hypothesis with an edge $\overline{a_1b_1}\subset \overline{ab}$ (Fig. 5.6). \\
\centerline{{\psfig{figure=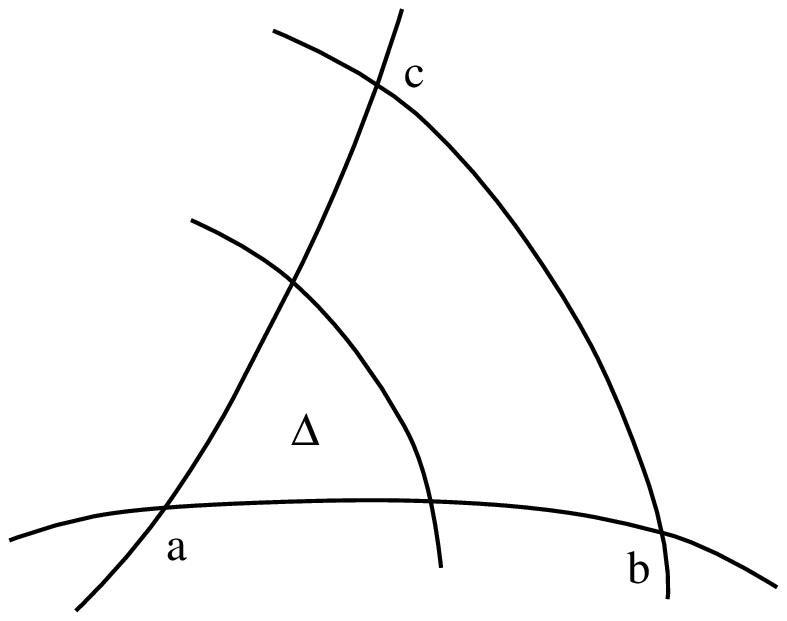,height=5.2cm}}}\ \\
\centerline{\footnotesize{Fig. 5.5}}\\
If $\gamma$ does not cut $\Delta_0$ or it cuts $\overline{a_1b_1}$ we are done (Fig 5.6). Therefore let us assume that $\gamma$ cuts $\overline{a_1c_1}$ 
(in $u_1$) and $\overline{b_1c_1}$ (in $w_1$). Let $\gamma$ cut $\overline{ab}$ in $u$ and $\overline{ac}$ in $w$ (Fig. 5.7). We have to consider two cases: \\
\indent (a) $\quad \overline{uu_1}\ \cap$ int $\Delta_0=\emptyset$ (so $\overline{ww_1}\ \cap$ int $\Delta_0=\emptyset$); Fig. 5.7.\\ 
\centerline{\epsfig{figure=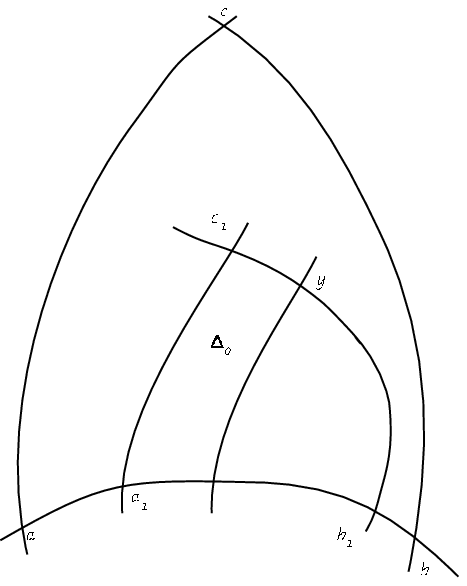, height=3.3in}}\\ 
\centerline{\footnotesize{Fig. 5.6}}\\ \ \\
\centerline{{\psfig{figure=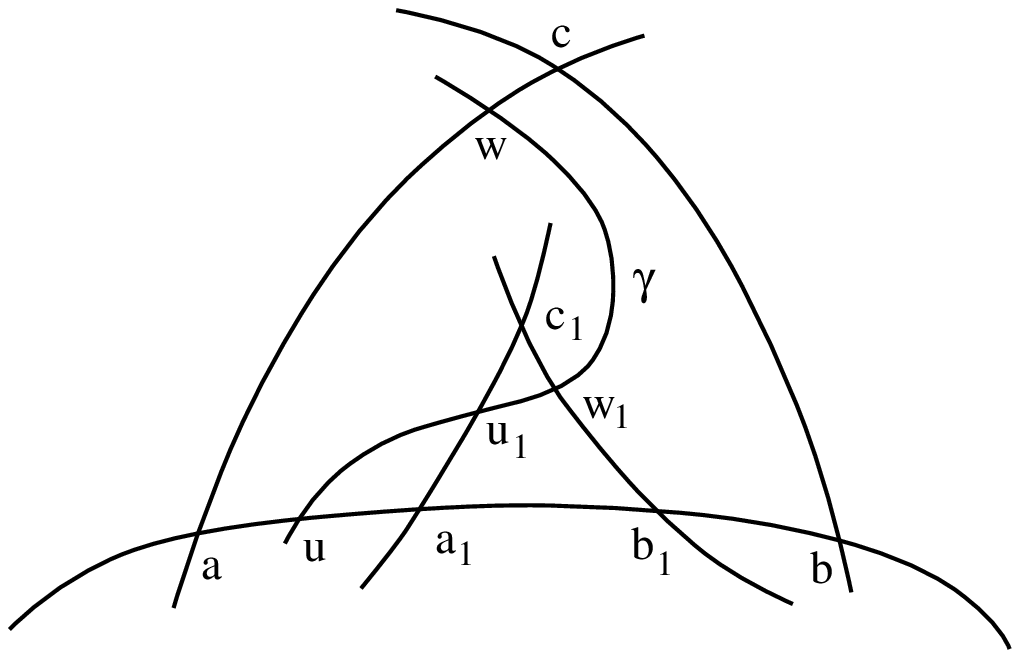,height=5.9cm}}}\ \\
\centerline{\footnotesize{Fig 5.7}}
Consider the 3-gon $ua_1u_1$. No arc can cut the edge $\overline{a_1u_1}$ so each arc which cuts the 3-gon $ua_1u_1$ cuts the edges $\overline{ua_1}$ 
and $\overline{uu_1}$. Furthermore this 3-gon is cut by less than $k+1$ arcs so by the inductive hypothesis there is a 3-gon $\Delta$ in $ua_1u_1$ 
with an edge on $\overline{ua_1}$ the interior of which is not cut by any arc. Then $\Delta$ satisfies the thesis of Proposition 5.2.\\
\indent (b) $\quad \overline{uw_1} \ \cap$ int$\Delta_0=\emptyset$ (so $\overline{wu_1}\ \cap$ int$\Delta_0 =\emptyset$). 
In this case we proceed like in case (a). \\
This completes the proof of Proposition 5.2 and hence the proof of Lemma 2.14.

\ \\
Department of Mathematics\\
Warsaw University\\
00-901 Warszawa, Poland\\ \ \\ \\
Added for e-print:\\
New address of J.~H.~Przytycki:\\
Department of Mathematics\\
The George Washington University\\
Washington, DC\\
{\tt przytyck@gwu.edu}\\
and University of Gda\'nsk

\end{document}